\documentclass[reqno, 11pt,  a4paper]{amsart}
\usepackage{color}
\usepackage{mathrsfs}
\newtheorem{thm}{Theorem}[section]

\theoremstyle{definition}

\theoremstyle{remark}
\newtheorem{rem}[thm]{Remark}
\newtheorem{example}[thm]{Example}
\numberwithin{equation}{section}
\newcommand{\Real}{\mathbb R}

\newcommand{\eps}{\varepsilon}
\newcommand{\diag}{\mathrm{diag}}

\newcommand{\F}{\mathscr{F}}
\newcommand{\s}{\mathbb{S}}

\newcommand{\one}[1]{\mathbf{1}_{\{#1\}}}
\newcommand{\bareone}[1]{\mathbf{1}_{#1}}

\renewcommand{\P}{\mathrm{P}}

\newcommand{\E}{\mathrm{E}}
\newcommand{\simplex}{\mathcal{S}^{d-1}}

\DeclareMathOperator*{\argmax}{\mathrm{argmax}}

\newcommand{\D}{\mathscr{D}}

\newcommand{\M}{\mathcal{M}}

\newcommand{\trace}{\mathrm{tr}}

\begin{document}

\title[]{On filtering of Markov chains in strong noise}%
\author{P.Chigansky}%
\address{Department of Mathematics, The Weizmann Institute of Science, Rehovot 76100, Israel}%
\email{pavel.chigansky@weizmann.ac.il}

\thanks{Research supported by a grant of the Israeli  Science Foundation}

\begin{abstract}
The filtering problem for finite state Markov chains is revisited in the low signal-to-noise regime.
We give a description of conditional measure concentration around the invariant distribution of the signal
and derive asymptotic expressions for the performance indices of the MMSE and MAP filtering estimates.

\end{abstract}
\maketitle
\section{Introduction}
Consider the discrete time signal/observation pair $(X,Y)=(X_n,Y_n)_{n\in\mathbb{Z}_+}$, where the {\em signal}
$X$ is a finite state Markov chain with values in a real alphabet $\s=\{a_1,...,a_d\}$, transition probabilities
$\lambda_{ij}=\P(X_n=a_j|X_{n-1}=a_i)$ and initial distribution $\nu$. The observation sequence $Y$ is generated by
\begin{equation}\label{Ygen}
Y_n = \sum_{i=1}^d \one{X_n=a_i}\xi_n(i), \quad n\ge 1
\end{equation}
where $\xi$ is a sequence of i.i.d. random vectors, independent of $X$. Without loss of generality the probability laws of the
entries of $\xi_1$ can be assumed to have densities $g_i(u)$, $i=1,...,d$, $u\in\Real$ with respect to a $\sigma$-finite measure $\psi(du)$ on $\Real$
(typically the Lebesgue measure or purely atomic measure). Hereafter all the random variables are assumed to be supported on a complete probability space $(\Omega,\F,\P)$.

This setting is often referred as Hidden Markov Model and is frequently encountered in information sciences
(see e.g. the recent survey \cite{EM}). An important statistical problem related to HMM is filtering, i.e. estimation of the
signal $X_n$, given the observation trajectory $Y$ up to time $n$. The main building blocks of this estimation problem are the conditional probabilities
$\pi_n(i)=\P(X_n=a_i|\F^Y_n)$, where $\F^Y_n=\sigma\{Y_m, m\le n\}$ is the $\sigma$-algebra of  events generated by the observations.
In particular the minimum mean square error (MMSE) and maximum a posterior probability (MAP) estimates of $X_n$ are given by
\begin{equation}\label{est}
\widehat{X}^{\mathrm{mse}}_n = \sum_{i=1}^d a_i \pi_n(i) \quad \text{and}\quad \widehat{X}^{\mathrm{map}}_n = \argmax_{a_i\in\s}\pi_n(i).
\end{equation}

The vector $\pi_n$ satisfies the recursive Bayes formula, called the {\em filtering} equation,
\begin{equation}
\label{filt}
\pi_n = \frac{G(Y_n)\Lambda^* \pi_{n-1}}{\big|G(Y_n)\Lambda^* \pi_{n-1}\big|}, \quad \pi_0=\nu,
\end{equation}
where $\Lambda^*$ is the transposed matrix of transition probabilities $\lambda_{ij}$, $G(y)$, $y\in\Real$ is the scalar matrix with
entries $g_i(y)$ and $|x|$ stands for the $\ell_1$-norm, i.e. $|x|=\sum_{i=1}^d|x_i|$. As usual we identify the probability measures and functions on $\s$ with vectors from the
simplex $\simplex=\{x\in\Real^d: x_i\ge 0, \sum_{i=1}^d x_i =1\}$ and  $\Real^d$ respectively and use the notation
$\eta(f)=\sum_{i=1}^d f(a_i)\eta_i=f^*\eta$ for $f:\s\mapsto\Real$ and $\eta\in\simplex$.

While the recursion \eqref{filt} provides an efficient way to calculate the estimates in \eqref{est}, no closed form formulae are known  for the
corresponding performance indices: the minimal mean square error
$$
\mathcal{E}_n= \min_{\theta\in \mathbb{L}^2(\Omega,\F^Y_n,\P)}\E\big(X_n-\theta\big)^2=\E\big(X_n-\widehat{X}^{\mathrm{mse}}_n\big)^2=
\E X_n^2 - \E \big(\pi_n(a)\big)^2
$$
and the minimum a posterior error probability
$$
\mathcal{P}_n = \min_{\theta\in\mathbb{L}^\infty(\Omega, \F^Y_n,\P)}\P\big(X_n\ne \theta\big) = 1-\E\max_{a_i\in \s}\pi_n(i).
$$
and hence  approximations of these quantities are of significant interest.

It is not hard to see that the random sequence $\pi_n$ is a Markov process with values in $\simplex$. Under mild assumptions it is also
a Feller process and hence it has at least one invariant measure $\M_\pi(d\eta)$ (on the Borel field of $\simplex$). The uniqueness of this measure
is not at all obvious and in fact may fail if no restrictions are imposed on the noise densities, even when the signal $X$ itself is ergodic
(see a discussion in \cite{BCL}).
Recall that a Markov chain on $\s$ is ergodic, if the limit probabilities $\mu_i:=\P(X_n=a_i)$, $i=1,...,d$ exist, are unique and positive.
The sufficient and necessary condition for ergodicity is that the matrix $\Lambda^q$ has positive entries for some integer $q\ge 1$ and
then $\mu$ is the unique solution of $\Lambda^*\mu=\mu$ in $\simplex$.
The invariant measure $\M_\pi$ of $\pi_n$ is unique, i.e. independent of $\nu$, if $X$ is ergodic and the noise densities are bounded and have the same
support (see \cite{Ch2}). In this case the limits
$$
\mathcal{E}:=\lim_{n\to\infty}\mathcal{E}_n \quad\text{and}\quad \mathcal{P} :=\lim_{n\to\infty}\mathcal{P}_n
$$
exist and do not depend on $\nu$.

Though these ``steady state'' optimal errors cannot be calculated exactly, they are amenable to asymptotic approximations, as the one
obtained by R.Khasminskii and O.Zeitouni in \cite{KZ} and G.Golubev in \cite{Gol}.
Suppose that the transition probabilities satisfy
$$
\lambda^\eps_{ij} = \begin{cases}
1-\eps\sum_{\ell\ne j}\lambda_{i\ell}, & i=j\\
\eps \lambda_{ij}, & i\ne j
\end{cases}
$$
with a small parameter $\eps\in(0,1)$, which controls the transitions rate of the corresponding {\em slow} chain $X^\eps_n$
(note that the invariant measure of $X^\eps$ does not depend on $\eps$ and equals $\mu$). The observation process $Y^\eps$ is given by \eqref{Ygen}
with $X$ replaced with $X^\eps$ and $\pi^\eps_n$ is the solution of \eqref{filt} with $Y$ and $\Lambda$ replaced with $Y^\eps$ and $\Lambda^\eps$
respectively.
It is proved in \cite{KZ}, that if all the Kullback-Leibler divergences
$$\D(g_i\parallel g_j)=
\int_\Real g_i(u)\log \frac{g_i}{g_j}(u)\psi(du)
$$
are finite and positive, the error probability\footnote{throughout
superscripts are added to various quantities to emphasize their
dependence on the corresponding parameter} $\mathcal{P}^\eps$
converges to zero as $\eps\to 0$ and
\begin{equation}
\label{KZas}
\mathcal{P}^\eps = \left(\sum_{i=1}^d \mu_i \sum_{j\ne i}\frac{\lambda_{ij}}{\D(g_j\parallel g_i)}\right)\eps\log\eps^{-1}\big(1+o(1)\big), \quad \eps\to 0.
\end{equation}
Similar asymptotic holds for the minimal mean square error as shown in \cite{Gol}:
\begin{equation}
\label{Gol}
\mathcal{E}^\eps = \left(\sum_{i=1}^d \mu_i\sum_{j\ne i}\frac{\lambda_{ij}}{\D(g_j\parallel g_i)}\big(a_i-a_j\big)^2\right)\eps\log\eps^{-1}\big(1+o(1)\big), \quad \eps\to 0.
\end{equation}
These results give an idea of how fast the invariant measure $\M_\pi^\eps(d\eta)$ concentrates around $\M_\pi^0(d\eta)=\sum_{i=1}^d\mu_i\delta_{p_i}(d\eta)$,
where $p_i$ are probability vectors with $1$ at the $i$-th entry.

In a sense the slow chain limit is the counterpart of the weak noise asymptotic $\sigma\to 0$ for the additive observation model
(cf. \eqref{Ygen})
\begin{equation}
\label{Y}
Y^\sigma_n = h(X_n) + \sigma\xi_n, \quad n\ge 1,
\end{equation}
where $\xi$ is a sequence of i.i.d. random variables, independent of $X$, $h$ is an $\s\mapsto\Real$ function and $\sigma$ is the constant, controlling
the noise intensity. Though less apparent in the discrete time setting, the analogy is complete for continuous time model, as explained in  Section \ref{sec2}
below. In this paper {\em the strong noise} asymptotic is addressed, when the filtering probabilities $\pi^\sigma_n$ converge to the
a priori distribution  of the signal $\nu_n=(\Lambda^*)^n\nu$ as $\sigma\to \infty$. Thus in the stationary case we deal with the
concentration of $\M_\pi^\sigma(d\eta)$ around $\M_\pi^\infty(d\eta)=\delta_\mu(d\eta)$ as $\sigma\to \infty$.
The precise formulation of the results is given in Section \ref{sec2}, which are proved in Sections \ref{sec3} and \ref{sec4}.

\section{Main results}\label{sec2}
\subsection{Discrete time} Let $(X,Y^\sigma)$ be the filtering model, with $X$ being a finite state Markov chain on $\s$ with transition probabilities
matrix $\Lambda$ and initial distribution $\nu$ and suppose that $Y^\sigma$ is generated by \eqref{Y}.
\begin{thm}\label{thm1}
Assume that the probability law of $\xi_1$ has a bounded twice
continuously differentiable density $g(u)$ with respect to the
Lebesgue measure on $\Real$ with bounded continuous derivatives.
Then the solution of \eqref{filt} converges to $\nu_n =
\big(\Lambda^{*}\big)^n\nu$ as $\sigma\to \infty$ and
$$
\sigma\big(\pi^\sigma_n-\nu_n\big)\xrightarrow[\sigma\to\infty]{\P-a.s.} Z_n,\quad n\ge 0
$$
where $Z_n$  satisfies
\begin{equation}
\label{Zeq}
Z_n = \Lambda^* Z_{n-1} - \big(\diag(\nu_n)-\nu_n\nu^*_n\big)h\frac{g'(\xi_n)}{g(\xi_n)}, \quad Z_0=0.
\end{equation}
\end{thm}
The following two theorems give asymptotic expressions for $\mathcal{E}^\sigma$ and $\mathcal{P}^\sigma$.
\begin{thm}\label{cor1}
Assume that $X$ is an ergodic chain and $g$ satisfies the following conditions
\begin{enumerate}
\renewcommand{\theenumi}{$\mathbf{a}_{\text{\arabic{enumi}}}$}
\item \label{a1} $g(u)$ does not vanish on $\Real$, is bounded and has two bounded derivatives
\item \label{a2} there is a $\delta>0$, so that
$$
\int_{-\infty}^\infty \left(\frac{g'(x)}{\min_{|u|\le \delta}g(x+u)}\right)^2g(x)dx<\infty,
$$
and
$$
\int_{-\infty}^\infty \left(\frac{\max_{|v|\le \delta}|g''(x+v)|}{\min_{|u|\le \delta}g(x+u)}\right)^2g(x)dx<\infty.
$$
\end{enumerate}
Let $I$ denote the Fisher information of $g$:
$$
I=\int_{-\infty}^\infty \frac{\big(g'(x)\big)^2}{g(x)}dx<\infty.
$$
Then the algebraic Lyapunov equation
\begin{equation}
\label{lyap}
P = \Lambda^* P\Lambda + \big(\diag(\mu)-\mu\mu^*\big)h I h^*\big(\diag(\mu)-\mu\mu^*\big)
\end{equation}
has a unique solution $P$ in the class of nonnegative definite matrices with $\sum_{i,j}P_{ij}=0$ and
\begin{equation}
\label{eqcor1}
\lim_{\sigma\to\infty}\sigma^2\big(\mathcal{E}^\infty-\mathcal{E}^\sigma \big) = a^*P a,
\end{equation}
where $a$ is a vector with entries $a_1,...,a_d$ and $\mathcal{E}^\infty=\mu(a^2)-\mu^2(a)$ is the a priori mean square error.
\end{thm}

\begin{rem}\label{rem}
The assumption \eqref{a1} and ergodicity of $X$ guarantee uniqueness of the invariant measure $\M_\pi^\sigma(d\eta)$ (see \cite{Ch2}).
The assumption \eqref{a2} is satisfied for many frequently encountered densities. For Gaussian density $g(x)=(2\pi)^{-1/2}\exp\{-x^2/2\}$
$$
\frac{|g'(x)|}{\min_{|u|\le \delta}g(x+u)} = \frac{|x|e^{-x^2/2}}{\min_{|u|\le \delta}e^{-(x+u)^2/2}}= \
\frac{|x|}{\min_{|u|\le \delta}e^{-xu-u^2/2}} \le
\frac{|x|}{e^{-|x|\delta-\delta^2/2}}
$$
and hence
$$
\int_{-\infty}^\infty \left(\frac{|g'(x)|}{\min_{|u|\le \delta}g(x+u)}\right)^pg(x)dx \le
\int_{-\infty}^\infty \left(\frac{|x|}{e^{-|x|\delta-\delta^2/2}}\right)^p e^{-x^2/2}dx<\infty
$$
for any $p\ge 0$ and not just $p=2$ as required by the first part of \eqref{a2}.
Similarly
$$
\frac{|g''(x+v)|}{g(x+u)} = \frac{\max_{|v|\le \delta}|(x+v)^2-1|e^{-(x+v)^2/2}}{\min_{|u|\le \delta} e^{-(x-u)^2/2}}\le
(2x^2+2\delta^2+1)e^{2|x|\delta+\delta^2}
$$
and the second condition of \eqref{a2} holds with any power $p\ge 0$ as well. It is not hard to verify that \eqref{a2} also holds for e.g.
Cauchy density $g(x)=\pi^{-1}(1+x^2)^{-1}$, which fails to have the first moment.

\end{rem}

\begin{thm}\label{cor2}
Assume that $X$ is ergodic and $\xi_1$ is a standard Gaussian random variable, then for any continuous function $F:\Real^d\mapsto\Real$,
growing not faster than polynomially,
\begin{equation}\label{Lplim}
\int_{\simplex}F\big(\sigma\big(\eta-\mu\big)\big)\M^\sigma_\pi(d\eta)\xrightarrow{\sigma\to \infty} \E F(Z),
\end{equation}
where $Z$ is a zero mean Gaussian vector with covariance matrix $P$, defined by \eqref{lyap} with $I\equiv 1$.
In particular
\begin{equation}\label{MAP}
\lim_{\sigma\to \infty}\sigma \big(\mathcal{P}^\infty - \mathcal{P}^\sigma\big) = \E\max_{j\in \mathcal{J}} Z_j,
\end{equation}
where $\mathcal{P}^\infty:=1-\max_{a_i\in\s}\mu_i$ is the a priori error probability and
$\mathcal{J}=\{i:\mu_i=\max_{j}\mu_j\}$. If $\mu$ has a unique maximal atom, then for any integer
$p\ge 1$
$$
\lim_{\sigma\to \infty}\sigma^p \big(\mathcal{P}^\infty - \mathcal{P}^\sigma\big)=0.
$$
\end{thm}
If the maximal atom of $\mu$ is not unique, then the right hand side of \eqref{MAP} is positive in general as
the following example demonstrates.
\begin{example}
Let $X$ be a binary chain with the transition matrix
$$
\Lambda=\begin{pmatrix}
\lambda &1-\lambda\\
1-\gamma & \gamma
\end{pmatrix}, \quad \lambda,\gamma\in (0,1).
$$
The equation \eqref{lyap} is one dimensional and $P:=P_{11}=P_{22}=-P_{12}=-P_{21}$ satisfies
$$
P = P(1-\lambda-\gamma)^2 + \mu_1^2\mu_2^2(h_1-h_2)^2
$$
and hence
$$
P=\frac{\mu_1^2\mu_2^2(h_1-h_2)^2}{(\lambda +\gamma)(1-\lambda+1-\gamma)}=
\frac{(1-\lambda)^2(1-\gamma)^2(h_1-h_2)^2}{(\lambda +\gamma)(1-\lambda+1-\gamma)^5}.
$$
Now by Theorem \ref{cor1},
$$
\lim_{\sigma\to\infty}\sigma\big(\mathcal{E}^\infty-\mathcal{E}^\sigma \big) = (a_1-a_2)^2 P.
$$
By Theorem \ref{cor2}, if $\gamma\ne \lambda$
$$
\lim_{\sigma\to \infty}\sigma^p \big(\mathcal{P}^\infty - \mathcal{P}^\sigma\big) = 0, \quad p\ge 1
$$
and if $\gamma=\lambda$,
\begin{multline*}
\lim_{\sigma\to \infty}\sigma \big(\mathcal{P}^\infty - \mathcal{P}^\sigma\big) =
\E \max(Z,-Z)=\E |Z| =\\ 2\sqrt{P} \int_0^\infty\frac{x}{\sqrt{2\pi}} e^{-x^2/2}dx
=\frac{|h_1-h_2|}{4\sqrt{\lambda(1-\lambda)}}\cdot 0.3839...
\end{multline*}
\end{example}
\subsection{Continuous time}
The continuous time analogue of the aforementioned setting consists of a time homogeneous Markov chain with values in $\s$, transition
intensities $\lambda_{ij}$ and initial distribution $\nu$ and the observation process $Y^\sigma=(Y^\sigma)_{t\in\Real_+}$ satisfying
$$
Y^\sigma_t = \int_0^t h(X_s)ds + \sigma B_t, \quad t\ge 0,
$$
where $h$ is an $\s\mapsto\Real$ function, $\sigma>0$ is a real constant and $B=(B_t)_{t\ge 0}$ is a Brownian motion, independent of $X$.
We treat the continuous time case separately and hence use the same notations for transition intensities and transition probabilities, etc.

The vector of conditional probabilities $\pi_t$ satisfies the Wonham filtering It\^o stochastic differential equation (\cite{W}, see also \cite{LS})
\begin{equation}\label{W}
d\pi_t = \Lambda^* \pi_t dt + \sigma^{-2} \big(\diag(\pi_t)-\pi_t\pi_t^*\big)h\big(dY_t-\pi_t(h)dt\big),
\end{equation}
subject to $\pi_0=\nu$, where $\Lambda$ is the transition intensities matrix, $\diag(x)$, $x\in\Real^d$ stands for the scalar matrix with $x_i$ on
the diagonal, $h$ is a column vector with entries $h(a_i)$ and $x^*$ is the transposed of $x$.

Recall that $X=(X_t)_{t\in\Real_+}$ is ergodic, if $\exp(\Lambda)$ has positive entries or equivalently if all of its states communicate.
For ergodic chains the Markov process $\pi_t$ has a unique invariant measure $\M^\sigma_\pi(d\eta)$ for any $\sigma>0$ (see \cite{Ch,Ch2}).
In the case $d=2$ the exact expressions are known for both $\mathcal{P}$ and $\mathcal{E}$ in terms of integrals with respect to
the density of $\M_\pi(d\eta)$, which can be explicitly found by solving the corresponding Kolmogorov-Fokker-Plank equation (see \cite{W}, \cite{LS}).
In higher dimension the closed form solution for  KFP equation is unavailable, which makes the direct analysis of \eqref{W} intractable.

The slow chain $X^\eps$ is obtained by the time scaling $X^\eps_t=X_{\eps t}$, $t\ge 0$ and its transition intensities matrix
equals $\eps\Lambda$. As in the discrete time the invariant measure of $X^\eps$ is independent of $\eps$ and solves $\Lambda^*\mu=0$ in $\simplex$.
The asymptotic expressions \eqref{KZas} and \eqref{Gol} remain valid with $\D(g_i\parallel g_j)$ replaced by $\big(h_i-h_j\big)^2/2$ (see \cite{KZ}, \cite{Gol}).
It is not hard to see, either by appropriate time change or directly from the KFP equation, that the weak noise asymptotic $\sigma\to 0$ is obtained by
replacing $\eps=\sigma^2$ in these expressions.
The strong noise asymptotic also turns to be similar to the discrete time case:
\begin{thm}\label{thm1ct}
The solution of \eqref{W} converges to $\nu_t=e^{\Lambda^* t}\nu$ as $\sigma\to\infty$ and for any $p\ge 1$
$$
\sigma\big(\pi^\sigma_t-\nu_t\big) \xrightarrow[\sigma\to\infty]{\mathbb{L}^p} Z_t, \quad t\ge 0,
$$
where $Z_t$ is the Gaussian diffusion process:
\begin{equation}\label{Zeqct}
dZ_t = \Lambda^* Z_t dt + \big(\diag(\nu_t)-\nu_t\nu_t^*\big)hd\bar{B}_t, \quad Z_0=0,
\end{equation}
with $\bar{B}=\sigma^{-1}\Big(Y^\sigma_t-\int_0^t \pi^\sigma_s(h)ds\Big)$ being the innovation Brownian motion.
If $X$ is ergodic, the algebraic Lyapunov equation
\begin{equation}
\label{alyap}
0=\Lambda^* P +P\Lambda + \big(\diag(\mu)-\mu \mu^*\big)hh^*\big(\diag(\mu)-\mu \mu^*\big)
\end{equation}
has a unique solution $P$ in the class of nonnegative definite matrices satisfying $\sum_{ij}P_{ij}=0$
and for any $F:\Real^d\mapsto\Real$, growing not faster than polynomially,
\begin{equation}
\label{mrs}
\int_{\simplex}F\big(\sigma\big(\eta-\mu\big)\big)\M^\sigma_\pi(d\eta)\xrightarrow{\sigma\to \infty} \E F(Z),
\end{equation}
where $Z$ is a zero mean Gaussian random vector with covariance matrix $P$.
\end{thm}
Theorems \ref{cor1} and \ref{cor2} remain valid in continuous time case with obvious adjustments, namely
$$
\lim_{\sigma\to\infty}\sigma^2\big(\mathcal{E}^\infty-\mathcal{E}^\sigma \big) = a^*P a
\quad \text{and}\quad
\lim_{\sigma\to \infty}\sigma \big(\mathcal{P}^\infty - \mathcal{P}^\sigma\big) = \E\max_{j\in \mathcal{J}} Z_j,
$$
where $P$ is the solution of \eqref{alyap} and $Z$ is the Gaussian vector defined in Theorem \ref{thm1ct}.
If the maximal atom of $\mu$ is unique,
$$
\lim_{\sigma\to \infty}\sigma^p \big(\mathcal{P}^\infty - \mathcal{P}^\sigma\big)=0,\quad p\ge 1.
$$

\section{Proofs in discrete time}\label{sec3}

\subsection{Proof of Theorem \ref{thm1}}

The entries of the diagonal matrix $G(y)$ in \eqref{filt} in case of the observations \eqref{Y} have the form
$$
g\left(\frac{y-h_i}{\sigma}\right), \quad i=1,...,d.
$$
To emphasize the dependence on $\sigma$ write $G^\sigma(y)$ and let
$$T^\sigma(y)x=G^\sigma(y)\Lambda^*x/|G^\sigma(y)\Lambda^*x|, \quad y\in \Real, \ x\in\simplex.$$
Since the density $g(u)$ is continuous, for any $i=1,...,d$
$$
g\left(\frac{Y^\sigma_n-h_i}{\sigma}\right) =
g\left(\xi_n+\frac{h(X_n)-h_i}{\sigma}\right) \xrightarrow[\sigma\to\infty]{\P-a.s.} g(\xi_n),
$$
and hence  $\lim_{\sigma\to \infty}T^\sigma(Y_n) x = \Lambda^* x$, $\P$-a.s. for any $x\in\simplex$. Then for any fixed $n\ge 1$
\begin{equation}\label{pipi}
\pi^\sigma_n = T^\sigma(Y^\sigma_n)\circ \cdots \circ T^\sigma(Y^\sigma_1)\circ \nu\xrightarrow[\sigma\to\infty]{\P-a.s.} (\Lambda^*)^n\nu = \nu_n.
\end{equation}
Since both $\pi^\sigma_n$ and $\nu_n$ are bounded, the convergence also holds in $\mathbb{L}^p$ for any $p\ge 1$.

Let $q^\sigma_n$ be the solution of
\begin{equation}
\label{qeq}
q^\sigma_n = \Lambda^* q^\sigma_{n-1} - \sigma^{-1} \big(\diag(\nu_n)-\nu_n\nu_n^*\big) h \frac{g'(\xi_n)}{g(\xi_n)}, \quad q^\sigma_0=\nu.
\end{equation}
The process $\Delta^\sigma_n = \sigma\big(\pi^\sigma_n-q^\sigma_n\big)$ satisfies
\begin{equation}\label{DD}
\Delta^\sigma_n = \Lambda^* \Delta^\sigma_{n-1} +
\sigma\bigg(\frac{G^\sigma(Y^\sigma_n)\Lambda^*\pi^\sigma_{n-1}}{|G^\sigma(Y^\sigma_n)\Lambda^*\pi^\sigma_{n-1}|}-\Lambda^* \pi^\sigma_{n-1}\bigg)
+\big(\diag(\nu_n)-\nu_n\nu_n^*\big) h \frac{g'(\xi_n)}{g(\xi_n)},
\end{equation}
subject to $\Delta^\sigma_0=0$.
Denote $\pi^\sigma_{n|n-1}=\Lambda^*\pi^\sigma_{n-1}$, then
\begin{align*}
& \frac{g\big(\sigma^{-1}(Y^\sigma_n-h_i)\big)\pi^\sigma_{n|n-1}(i)}
{\big|G^\sigma(Y^\sigma_n)\pi^\sigma_{n|n-1}\big|}-\pi^\sigma_{n|n-1}(i)=\\
& \frac{
g\big(\xi_n + \sigma^{-1}(h(X_n)-h_i) \big)-
\sum_{j=1}^d g\big(\xi_n + \sigma^{-1}(h(X_n)-h_j)\big)\pi^\sigma_{n|n-1}(j)
}
{\sum_{j=1}^d g\big(\xi_n + \sigma^{-1}(h(X_n)-h_j)\big)\pi^\sigma_{n|n-1}(j)}\pi^\sigma_{n|n-1}(i) =\\
& \frac{
-\sigma^{-1}g'(\xi_n)\Big(h_i-\sum_{j=1}^d h_j\pi^\sigma_{n|n-1}(j)\Big) + K_n \sigma^{-2}
}
{\sum_{j=1}^d g\big(\xi_n + \sigma^{-1}(h(X_n)-h_j)\big)\pi^\sigma_{n|n-1}(j)}\pi^\sigma_{n|n-1}(i),
\end{align*}
where $K_n$ are  bounded random variables  (recall that $g''(u)$ is assumed bounded).
Hence by \eqref{pipi} and continuity of $g$
$$
\sigma\left(\frac{G^\sigma(Y^\sigma_n)\pi^\sigma_{n|n-1}}{|G^\sigma(Y^\sigma_n)\pi^\sigma_{n|n-1}|}-\pi^\sigma_{n|n-1}\right)\xrightarrow[\sigma\to\infty]{\P-a.s.}
-\big(\diag(\nu_n)-\nu_n\nu_n^*\big) h \frac{g'(\xi_n)}{g(\xi_n)}.
$$
Iterating  \eqref{DD} one gets
$
\lim_{\sigma\to \infty}\Delta^\sigma_n=0
$,
$ \P-a.s.$ for any fixed $n\ge 0$ and the statement of the theorem follows:
$$
\sigma\big(\pi^\sigma_n-\nu_n\big) = \sigma\big(\pi^\sigma_n-q^\sigma_n\big) + \sigma\big(q^\sigma_n-\nu_n\big)\xrightarrow[\sigma\to \infty]{\P-a.s.} Z_n,
$$
where $Z_n:=\sigma(q^\sigma_n-\nu_n)$ clearly satisfies \eqref{Zeq}, which doesn't depend on $\sigma$.\qed

\subsection{Proof of Theorem \ref{cor1}}
Note that to verify \eqref{eqcor1} one should first take the limit $n\to\infty$ and then $\sigma\to\infty$ and thus
cannot use the statement of Theorem \ref{thm1} {\em per se}. The proof relies on stability property of the matrix $\Lambda$, provided by
ergodicity of $X$.

As shown in \cite{Ch2} the Markov process $(X,\pi^\sigma)$ has the unique invariant measure $\M^\sigma(dx,d\eta)$ if $X$ is ergodic and
assumption \eqref{a1} is satisfied. In particular  $\M^\sigma_\pi(d\eta)= \sum_{i=1}^d \M^\sigma(\{a_i\}, d\eta)$. If the equation \eqref{filt} and $X$ is started from a random variable with
distribution $\M^\sigma(dx,d\eta)$, the process $\pi^\sigma = (\pi^\sigma_n)_{n\ge 0}$ is stationary, which is assumed hereafter.

As in \eqref{pipi}
$
\lim_{\sigma\to\infty}\pi^\sigma_n  = (\Lambda^*)^n\pi^\sigma_0,\ \P-a.s.
$
and so for any $\eps>0$ and any $m\ge 0$
\begin{multline}\label{statlim}
\lim_{\sigma\to\infty}\P\big(|\pi^\sigma_0-\mu|\ge \eps \big)=\lim_{\sigma\to\infty}\P\big(|\pi^\sigma_m-\mu|\ge \eps\big)\le \\
\lim_{\sigma\to\infty}\P\big(|\pi^\sigma_m-(\Lambda^*)^m\pi^\sigma_0|\ge \eps/2\big) +\varlimsup_{\sigma\to\infty}\P\big(|(\Lambda^*)^m\pi^\sigma_0-\mu|\ge \eps/2\big)
\xrightarrow{m\to\infty} 0
\end{multline}
where the latter convergence holds since $(\Lambda^*)^nx\to \mu$ for all $x\in\simplex$ by ergodicity of $X$.

Let $q^\sigma_n$ denote the solution of (cf. \eqref{qeq})
$$
q^\sigma_n=\Lambda^* q^\sigma_{n-1} - \sigma^{-1} \big(\diag(\mu)-\mu\mu^*\big)h\frac{g'(\xi_n)}{g(\xi_n)}, \quad q^\sigma_0=\mu.
$$
We use the notations, introduced in the previous section, to denote random processes, playing the same role as in the proof
of Theorem \ref{thm1}, but defined differently to fit the stationary setup under consideration.

The process $\Delta^\sigma_n = \sigma(\pi^\sigma_n-q^\sigma_n)$ satisfies (cf. \eqref{DD})
\begin{multline*}
\Delta^\sigma_n = \Lambda^* \Delta^\sigma_{n-1} +
\sigma\bigg(\frac{G^\sigma(Y^\sigma_n)\Lambda^*\pi^\sigma_{n-1}}{|G^\sigma(Y^\sigma_n)\Lambda^*\pi^\sigma_{n-1}|}-\Lambda^* \pi^\sigma_{n-1}\bigg)
+ \\
\big(\diag(\mu)-\mu\mu^*\big) h \frac{g'(\xi_n)}{g(\xi_n)} := \Lambda^* \Delta^\sigma_{n-1}  + \theta^\sigma_n
\end{multline*}
subject to $\Delta^\sigma_0=\sigma(\pi^\sigma_0-\mu)$. Note that since the Fisher information is finite and $\pi^\sigma_n$ is stationary,
for any fixed $\sigma>0$,
$
\E \theta^\sigma_n \theta^{\sigma*}_n = \E \theta^\sigma_0 \theta^{\sigma*}_0 := \Gamma^\sigma
$
and hence $Q^\sigma_n:=\E \Delta^\sigma_n \Delta^{\sigma*}_n$ satisfies
\begin{equation}
\label{Qlyap}
Q^\sigma_n = \Lambda^* Q^\sigma_{n-1} \Lambda + \Gamma^\sigma, \quad n\ge 1,
\end{equation}
subject to $Q^\sigma_0=\sigma^{2}\E(\pi^\sigma_0-\mu)(\pi^\sigma_0-\mu)^*$. If $X$ is ergodic, $\Lambda^*$ is a stability matrix, when
restricted to the subspace $\{x\in\Real^d:\sum_{i=1}^d x_i=0\}$. Since $\Delta^\sigma_n$ belongs to this subspace for all $n\ge 0$, the Lyapunov equation
\eqref{Qlyap} has a bounded solution, which converges to the unique limit
$
Q^\sigma = \sum_{m=0}^\infty \Lambda^{*m}\Gamma^\sigma \Lambda^m.
$

For brevity set $\pi^\sigma_{1|0}=\Lambda^*\pi^\sigma_0$ and define
$$
a^\sigma :=
\sigma\bigg(\frac{G^\sigma(Y^\sigma_1)\pi^\sigma_{1|0}}{\big|G^\sigma(Y^\sigma_1)\pi^\sigma_{1|0}\big|}-\pi^\sigma_{1|0}\bigg).
$$
Then
\begin{align*}
& a^\sigma(i) =
 \sigma\left(\frac{g\big(\sigma^{-1}(Y^\sigma_1-h_i)\big)\pi^\sigma_{1|0}(i)}
{\big|G^\sigma(Y^\sigma_1)\pi^\sigma_{1|0}\big|}-\pi^\sigma_{1|0}(i)\right)=\\
& \sigma \frac{
g\big(\xi_1 + \sigma^{-1}(h(X_1)-h_i) \big)-
\sum_{j=1}^d g\big(\xi_1 + \sigma^{-1}(h(X_1)-h_j)\big)\pi^\sigma_{1|0}(j)
}
{\sum_{j=1}^d g\big(\xi_1 + \sigma^{-1}(h(X_1)-h_j)\big)\pi^\sigma_{1|0}(j)}\pi^\sigma_{1|0}(i) =\\
&
\frac{-g'(\xi_1)\big(h_i-\sum_{j=1}^d h_j\pi^\sigma_{1|0}(j)\big)}
{\sum_{j=1}^d g\big(\xi_1 + \sigma^{-1}(h(X_1)-h_j)\big)\pi^\sigma_{1|0}(j)}\pi^\sigma_{1|0}(i)
+ \\
& \frac{\pi^\sigma_{1|0}(i)}{2\sigma}\frac{
g''(\xi_1+ \alpha_i/\sigma)(h(X_1)-h_i)^2\beta_i-
\sum_{j=1}^dg''(\xi_1+ \alpha_j/\sigma)(h(X_1)-h_j)^2\beta_j\pi^\sigma_{1|0}(j)
}
{\sum_{j=1}^d g\big(\xi_1 + \sigma^{-1}(h(X_1)-h_j)\big)\pi^\sigma_{1|0}(j)}
\end{align*}
where the latter holds by the mean value theorem with $|\alpha_j|\le |h(X_1)-h_j|$ and $\beta_j\in[0,1]$.
Since $g''$ is bounded  and by \eqref{statlim} $\pi^\sigma_{1|0}\to \mu$ in probability as $\sigma\to\infty$
\begin{equation}\label{inprob}
a^\sigma(i)
\xrightarrow[\sigma\to\infty]{\P}
-\frac{g'(\xi_1)}{g(\xi_1)}\big(h_i-\sum_{j=1}^d h_j\mu_j\big)
\mu_i.
\end{equation}
Note that for $\sigma > \max_{i,j}|h_i-h_j|/\delta$,
$$
\frac{|g'(\xi_1)|}
{\sum_{j=1}^d g\big(\xi_1 + \sigma^{-1}(h(X_1)-h_j)\big)\pi^\sigma_{1|0}(j)} \le
\frac{|g'(\xi_1)|}
{\min_{|u|\le \delta}g\big(\xi_1 + u\big)},
$$
where by the assumption \eqref{a2} the right hand side is  square integrable.
Analogously for sufficiently small $\sigma$,
$$
\frac{|g''(\xi_1+ \alpha_i/\sigma)|}
{\sum_{j=1}^d g\big(\xi_1 + \sigma^{-1}(h(X_1)-h_j)\big)\pi^\sigma_{1|0}(j)}
\le
 \frac{\max_{|v|\le \delta}\big|g''(\xi_1+ v)\big|}
{
\min_{|u|\le \delta} g\big(\xi_1 + u\big)
},
$$
with a square integrable right hand side. Hence by the Lebesgue dominated convergence
\eqref{inprob} implies
$$
a^\sigma \xrightarrow[\sigma\to \infty]{\mathbb{L}^2}
-\big(\diag(\mu)-\mu\mu^*\big) h \frac{g'(\xi_1)}{g(\xi_1)}
$$
and in turn
\begin{equation}\label{Qsig}
\lim_{\sigma\to\infty}\trace(\Gamma^\sigma) =0\quad \implies \quad \lim_{\sigma\to\infty}Q^\sigma=0.
\end{equation}
On the other hand, the sequence $Z_n=\sigma(q^\sigma_n-\mu)$ does not depend on $\sigma$ and satisfies
\begin{equation}
\label{rhoeq}
Z_n=\Lambda^* Z_{n-1} - \big(\diag(\mu)-\mu\mu^*\big)h\frac{g'(\xi_n)}{g(\xi_n)}, \quad Z_0=0.
\end{equation}
Again by the stability property of  $\Lambda^*$ on the subspace $\{x\in\Real^d: \sum_{j=1}^d x_j =0\}$
$$
\lim_{n\to\infty}\E Z_nZ_n^* = P,
$$
where $P$ uniquely solves \eqref{alyap} in the class of nonnegative matrices with $\sum_{ij}P_{ij}=0$,
which is a well known property of the Lyapunov equation for stable matrices (see e.g. \cite{Bellman}).
Hence,
\begin{multline}\label{L2}
\sigma^2 \E(\pi^\sigma_0-\mu)(\pi^\sigma_0-\mu)^* - P =
\sigma^2 \E(\pi^\sigma_n-\mu)(\pi^\sigma_n-\mu)^* - \lim_{n\to\infty}Z_nZ_n^* = \\
\lim_{n\to\infty}\E \Delta^\sigma_n \Delta^{\sigma *}_n + \lim_{n\to\infty}\E (\Delta^\sigma_n Z^*_n+Z_n^* \Delta^\sigma_n)
\xrightarrow{\sigma\to \infty} 0
\end{multline}
where the latter convergence holds by \eqref{Qsig}.
This in turn implies \eqref{eqcor1}:
$$
\sigma^2\big(\mathcal{E}^\infty-\mathcal{E}^\sigma\big)=
\sigma^2\E a^*(\pi^\sigma_0-\mu)(\pi^\sigma_0-\mu)^*a
\xrightarrow{\sigma\to\infty}a^*Pa.\qed
$$

\subsection{Proof of Theorem \ref{cor2}}
For standard Gaussian $\xi_1$, $p'(x)/p(x)=-x$ and $I=1$, hence the process $Z_n$, defined in \eqref{rhoeq} is Gaussian.
By the Remark \ref{rem} Gaussian density satisfies the assumption \eqref{a2} of Theorem \ref{cor1} with square integrability
replaced by integrability to any power $p\ge 1$ and hence, similarly to \eqref{L2},
for any continuous $F:\Real^d\mapsto\Real$ with the norm bounded by a polynomial function of any finite order
$$
\E F\big(\sigma(\pi^\sigma_0-\mu)\big)\xrightarrow{\sigma\to\infty} \E F(Z),
$$
where $Z$ is a zero mean Gaussian vector with covariance matrix $P$, defined by \eqref{lyap}.

Let $\mathcal{J}=\{i:\mu_i=\max_j\mu_j\}$ and assume $\mu_1\in\mathcal{J}$ for definiteness.
Then
\begin{multline*}
\sigma \big(\mathcal{P}^\infty - \mathcal{P}^\sigma\big)=\sigma\big(\E \max_{a_i\in\s}\pi^\sigma_0(i)-\max_{i}\mu_i\big) =
\E \sigma\max_{a_i\in\s}(\pi^\sigma_0(i)-\mu_1) =\\ \E \max_{a_i\in\s}\big(\sigma(\pi^\sigma_0(i)-\mu_i)+\sigma(\mu_i-\mu_1)\big)
\xrightarrow{\sigma\to\infty}\E \max_{j\in \mathcal{J}}Z_j,
\end{multline*}
where the convergence holds by \eqref{Lplim}, since $\max_i(x_i)$, $x\in\Real^d$ is a continuous function and $\mu_i-\mu_1<0$ for $i\not\in\mathcal{J}$.

Suppose now that $\mu_1$ is the unique maximal atom of $\mu$ and let $r = \min_{j\ne 1}|\mu_1-\mu_j|>0$. Let $A_\sigma:=\{|\pi^\sigma_0-\mu|\le r/2\}$,
$\bareone{A_\sigma}$ be the indicator function of $A_\sigma$ and $A^c_\sigma=\Omega\backslash A_\sigma$.
Then
$$
\max_{a_i\in\s}\pi^\sigma_0(i) = \bareone{A_\sigma}\pi^\sigma_0(1) + \bareone{A^c_\sigma}\max_{a_i\in\s}\pi^\sigma_0(i)=
\pi^\sigma_0(1)+ \bareone{A^c_\sigma}\big(\max_{a_i\in\s}\pi^\sigma_0(i)-\pi^\sigma_0(1)\big),
$$
Hence for any two  integers $q>p\ge 1$,
\begin{multline*}
\sigma^{p}\big|\E \max_{a_i\in\s}\pi^\sigma_0(i)-\max_{i}\mu_i\big| =\\
\sigma^{p}\big|\big(\E \pi^\sigma_0(1)-\mu_1\big) +
\E \bareone{A^c_\sigma}\big(\max_{a_i\in\s}\pi^\sigma_0(i)-\pi^\sigma_0(1)\big)\big| =\\
\sigma^p\E \bareone{A^c_\sigma}\big|\max_{a_i\in\s}\pi^\sigma_0(i)-\pi^\sigma_0(1)\big|\le
2\sigma^p \frac{\E |\pi^\sigma_0-\mu|^q }{(r/2)^q}
\xrightarrow{\sigma\to \infty} 0,
\end{multline*}
since by \eqref{Lplim}, the limit
$\lim_{\sigma\to\infty}\sigma^q \E |\pi^\sigma_0-\mu|^q $ exists and is finite.\qed
\section{Proofs in continuous time}\label{sec4}

The continuous time filter \eqref{W} has a convenient innovation structure, which  simplifies the proofs, which are only sketched below.
\subsection{Proof of Theorem \ref{thm1ct}}
Since $\nu_t=\exp(\Lambda^*t)\nu$ solves $\dot{\nu}_t=\Lambda^* \nu_t$, $\nu_0=\nu$,
the process $\delta^\sigma_t:=\pi^\sigma_t-\nu_t$ satisfies
$$
d\delta^\sigma_t = \Lambda^* \delta^\sigma_t dt + \sigma^{-1}\big(\diag(\pi^\sigma_t)-\pi^\sigma_t\pi^{\sigma *}_t\big)hd\bar{B}_t, \quad \delta^\sigma_0=0,
$$
and hence
$$
\delta^\sigma_t = \sigma^{-1}\int_0^t e^{\Lambda^*(t-s)}\big(\diag(\pi^\sigma_s)-\pi^\sigma_s\pi^{\sigma *}_s\big)hd\bar{B}_s.
$$
Since the integrand is continuous and bounded for any $t\ge 0$,
$$
\lim_{\sigma\to\infty}\delta^\sigma_t=0, \quad \P-a.s.
$$
The convergence holds in $\mathbb{L}^p$, $p\ge 1$  as well, since the integrand of the stochastic integral is uniformly bounded in $\sigma$ and hence
$|\delta^\sigma_t|$ is uniformly integrable to any power as $\sigma\to \infty$.
Let $q^\sigma_t$ be solution of the linear SDE
$$
dq^\sigma_t = \Lambda^* q^\sigma_t dt + \sigma^{-1}\big(\diag(\nu_t)-\nu_t\nu_t^*\big)hd\bar{B}_t, \quad q^\sigma_0=\nu.
$$
The process $\Delta^\sigma_t=\sigma(\pi^\sigma_t-q^\sigma_t)$ satisfies
\begin{equation}\label{Deltact}
d\Delta^\sigma_t = \Lambda^* \Delta^\sigma_tdt + \big(\Gamma(\pi^\sigma_t)-\Gamma(\nu_t)\big)hd\bar{B}_t, \quad \Delta^\sigma_0=0,
\end{equation}
where $\Gamma(x)=\diag(x)-xx^*$ is defined for brevity. Then
$$
\Delta^\sigma_t = \int_0^t e^{\Lambda^*(t-s)}\big(\Gamma(\pi^\sigma_s)-\Gamma(\nu_s)\big)hd\bar{B}_s\xrightarrow{\sigma\to\infty}0, \quad \P-a.s\ \text{and in\ }\mathbb{L}^p,
$$
since $\Gamma(\cdot)$ is continuous, $\pi^\sigma_t$ and $\nu_t$ are bounded and $\pi^\sigma_t\to\nu_t$ $\P$-a.s. as $\sigma\to\infty$.
The process  $Z_t=\sigma(q^\sigma_t-\nu_t)$ satisfies \eqref{Zeqct} and thus
$$
\sigma\big(\pi^\sigma_t-\nu_t\big)=\sigma\big(\pi^\sigma_t-q^\sigma_t\big)+\sigma\big(q^\sigma_t-\nu_t\big)\xrightarrow[\sigma\to \infty]{\mathbb{L}^p}Z_t.
$$
If $X$ is ergodic, $\Lambda^*$ is a stability matrix on $\{x\in\Real: \sum_i x_i=0\}$, which is an invariant subspace of \eqref{Deltact} and \eqref{Zeqct}.
Then by the very same arguments, used in the proof of Theorem \ref{cor1}, and taking into account the integrability properties of the
stochastic integral with respect to $\bar{B}$, one verifies \eqref{mrs}. \qed

\subsection*{Acknowledgment} The author is grateful to Ofer Zeitouni and Tsachy Weissman for their useful comments about this paper.


\end{document}